\documentclass[10pt,twoside]{article}
\usepackage{graphicx}
\usepackage{amsmath}
\usepackage{Latex-document}

\newtheorem{theorem}{Theorem}[section]

\newcommand{\F}{{\bf F}}

\newcommand{\Tr}{\hbox{\rm Tr}}
\newcommand{\Norm}{\hbox{\rm Norm}}

\renewcommand{\thesection}{\arabic{section}}
\renewcommand{\thesubsection}{\thesection.\arabic{subsection}.}
\def\addsec{\addtocounter{section}{1} \setcounter{theorem}{0}}

\title{\bf Combinatorial Problems in Finite \vskip -2mm Geometry  and Lacunary Polynomials\vskip 6mm}

\author{{Aart Blokhuis}\thanks{Eindhoven University of Technology, P.O.\ Box 513, 5600 MB,
Eindhoven, Netherlands. E-mail: aartb@win.tue.nl}\vspace*{-0.5cm}}

\date{\vspace{-8mm}}

\begin{document}

\maketitle

\thispagestyle{first} \setcounter{page}{537}

\begin{abstract}

\vskip 3mm

We describe some combinatorial problems in finite projective planes and indicate how R\'edei's theory of lacunary
polynomials can be applied to them.

\vskip 4.5mm


\end{abstract}

\vskip 12mm

\section*{1. Introduction} \addsec

\vskip-5mm \hspace{5mm}

In 1991 I wrote a survey paper called Extremal Problems in Finite
Geometries \cite{Bl-ext}. It concerns among others problems of the
following type:

Given a set $B$ of points in a finite projective plane $\Pi$, with
the property that there is a restricted number of possibilities
for the size of the intersection of a line with $B$. What can be
concluded about the size and the structure of $B$.

The archetypal result is Segre's theorem \cite{Segre}:

{\it If \  $\Pi=PG(2,q)$, $q$ odd, and $B$ has at most two points
on a line, then $|B|\le q+1$ {\rm(this part is easy)}, and in case
of equality $B$ consists of the points of a conic.}

The problem becomes much more difficult when larger intersections
are allowed. A subset $B$ of $PG(2,q)$ of size $k$ having at most
$n$ points on a line, is called a $(k,n)$-arc.

A simple counting argument going back to Barlotti \cite{Barlotti}
gives that if $B$ is a $(k,n)$-arc,
 then
$k\le (n-1)(q+1)+1=nq-q+n$, and equality implies that $n|q$ and
all lines intersect $B$ in 0 or $n$ points.

An arc $B$ meeting the above upper bound is called a {\em maximal
arc}. The first non trivial case is $n=2$, and $q$ is even. In
this case $|B|=q+2$ is possible and the maximal arc is called a
hyperoval. In fact every $(q+1)$-arc can be extended to a
hyperoval by adding one point and classifying or trying to find
new hyperovals is one of the very active areas in finite geometry.
For a survey on the current situation we refer to the two papers
by Hirschfeld and Storme \cite{HiSt1,HiSt2}.

For $n=3$ non-existence of maximal arcs in $PG(2,9)$ was shown by
Cossu \cite{Cossu}, and later by Thas for general $q$
\cite{Thas1}. Very little is known however in this case. In fact,
if $k(q)$ stands for the maximal size of a $(k,3)$-arc in
$PG(2,q)$ then it is unknown whether limsup $k(q)/q>2$ or liminf
$k(q)<3$.

For $n=4$, and in fact for all $n=2^a$ and $q=2^b$, $b>a$ examples
are known, most of them due to Denniston \cite{Denniston} and Thas
\cite{Thas2,Thas3}.

For odd $n$ it was conjectured in \cite{Thas1} that maximal arcs
don't exist. This was finally proved by Ball, Blokhuis and
Mazzocca in 1996 \cite{BaBlMa}. A simplified proof appeared a year
later \cite{BaBl}.

 Now we turn to the case where $B$ intersects
every line in at least one point. Such a set is called a {\em
blocking set}. It is called non-trivial if it does not contain a
line. Here the classical result is due to Bruen \cite{Bruen}:

{\it A non-trivial blocking set $B$ in a projective plane of order
$q$ has size at least $q+\sqrt q+1$, with equality if and only if
$q$ is a square, and $B$ is a Baer subplane.}

 Our understanding
of the situation when $q$ is not a square has increased
dramatically in the last 10 years, from knowing very little to
more or less complete knowledge.

 Starting point was the
unexpectedly simple proof (in 1993) \cite{Blo94} that a
non-trivial blocking set in $PG(2,p)$, $p$ prime, has size at
least $3(p+1)/2$.

 The proof is based on properties of a certain
kind of lacunary polynomial (as introduced and studied by R\'edei
in \cite{Redei}) associated to the blocking set. The importance of
R\'edei's work in this area was realized soon after the appearance
of his book in 1970 notably in papers by Bruen and Thas,
\cite{BrTh}, where his result on the number of directions
determined by the graph of a function on a finite field was used.

 In contrast to the case that $B$ is an arc, we can still say
something if we require that $B$ intersects every line of the
plane in at least $t$ points, for some $t>1$. If $q$ is a square
then a natural candidate for $B$ is the union of $t$ disjoint Baer
subplanes. These can be found for all appropriate values of $t$
because it is possible to partition $PG(2,q)$ in $q-\sqrt q+1$
disjoint Baer subplanes.

Building on earlier work by Ball \cite{Ball96}, G\'acs, Sz\H{o}nyi
\cite{GaSz} and others Blokhuis, Storme and Sz\H{o}nyi showed that
for $t<q^{1/6}$ a $t$-fold blocking set in $PG(2,q)$ has size at
least $t(q+\sqrt{q}+1)$ with equality if and only if $B$ is the
union of $t$ disjoint Baer subplanes \cite{BlStSz}.

\section*{2. Directions} \addsec

\vskip-5mm \hspace{5mm}

Let $f\,:\,GF(q)\to GF(q)$. Define the set of directions
determined by $f$ to be
\[
D_f:=\left\{ {f(u)-f(v)\over u-v}\,|\,u,v\in GF(q), u\ne
v\right\}.
\]
We are interested in functions $f$ for which the set $D_f$ is
small. If $f$ is linear, then $D_f$ just consists of the slope of
the line defined by the graph of $f$. Our starting point will be
the following important result of R\'edei \cite{Redei}, p.\ 237,
Satz 24.
\begin{theorem}{\rm [R\'edei, 1970]}
Let $f:GF(q)\to GF(q)$ be a nonlinear function, where $q=p^n$, $p$
prime. Then $|D_f|$ is contained in one of the intervals
\[
\left(1+{q-1\over p^e+1},{q-1\over p^e-1}\right), e=1,\dots,[n/2];
\left( {q+1\over2},q\right).
\]
\end{theorem}

Examples of functions determining relatively few directions are
given by:
\begin{enumerate}
\item $f(x)=x^{q+1\over2}$, $q$ odd, $|D_f|=(q+3)/2$;
\item $f(x)=x^{p^e}$, where $e|n$, $|D_f|=(q-1)/( p^e-1)$;
\item $f(x)=\Tr(GF(q)\to GF(p^e))(x)$, $|D_f|=(q/ p^e)+1$.
\end{enumerate}

In all the examples $|D_f|$ is contained in one of the R\'edei
intervals corresponding to a subfield of $GF(q)$, i.e., $e|n$, and
the obvious question was whether this could be proved. In
\cite{BBS92} Blokhuis, Brouwer and Sz\H{o}nyi slightly improved
R\'edei's result but the real progress came with the paper
\cite{BBBSS} where not only it was proved that only the intervals
with $e|n$ occur, but the functions for which $|D_f|<(q+3)/2$
where essentially characterized:

\begin{theorem}{\rm [Ball, Blokhuis, Brouwer, Storme, Sz\H{o}nyi, 1999]}
Let $f:GF(q)\to GF(q)$, where $q=p^h$, $p$ prime, $f(0)=0$. Let
$N=|D_f|$. Let $e$ (with $(0\le e\le n)$ be the largest integer
such that each line with slope in $D_f$ meets the graph of $f$ in
a multiple of $p^e$ points. Then we have one of the following:
\begin{enumerate}
\item $e=0$ and $(q+3)/2\le N\le q+1$,
\item $e=1$, $p=2$ and $(q+5)/3\le N\le q-1$,
\item $p^e>2$, $e\,|\,n$, and $q/p^e+1\le N\le (q-1)/(p^e-1)$,
\item $e=n$ and $N=1$.
\end{enumerate}
Moreover, if $p^e>3$ or ($p^e=3$ and $N=q/3+1$), then $f$ is a
linear map on $GF(q)$ viewed as a vector space over $GF(p^e)$.
\end{theorem}

Very recently this result has been perfected by Simeon Ball,
removing the condition $p^e>2$ in the third case (and thus getting
rid of the second).

When we consider the set $B$ formed by the  $q$ points of the
graph of $f$ together with the $N=|D_f|$ points on the line at
infinity corresponding to the directions determined by $f$, we get
a blocking set. For if $l$ has a slope determined by $f$ then the
infinite point of $l$ belongs to $B$, and if not, then $l$ and its
parallels all contain precisely one point of the graph of $f$.

Conversely, if $B$ is a blocking set in $PG(2,q)$ of size $q+N$,
and there is a line intersecting $B$ in $N$ points, then it arises
from this construction. The blocking set is then called {\em of
R\'edei type}.

As mentioned in the introduction, the smallest non-trivial
blocking sets were characterized by Bruen \cite{Bruen} to be Baer
subplanes. They are of R\'edei type and correspond to the function
$x\mapsto x^{\sqrt q}$.

If the blocking set $B$ is of R\'edei type, then as a consequence
of the direction theorem above the structure is very special if
$N<(q+3)/2$, or equivalently if $|B|\le {3\over2}(q+1)$. An
important step towards showing that this is true in general is the
following result for planes of prime order already mentioned in
the introduction \cite{Blo94}.
\begin{theorem}{\rm [Blokhuis, 1994]}
Let $B$ be a blocking set in $PG(2,p)$, $p$ prime, not containing
a line. Then
\[
|B|\ge {3\over2}(p+1).
\]
\end{theorem}

The proof is based on properties of lacunary polynomials,
introduced and studied by R\'edei in \cite{Redei}. In the same
paper it is proved that a blocking set in $PG(2,p^3)$ has size at
least $p^3+p^2+1$. Recently Polverino \cite{P} has shown that
small blocking sets in $PG(2,p^3)$ are all of R\'edei type, and
the possible sizes are $p^3+p^2+1$ and $p^3+p^2+p+1$
(corresponding to examples 2 and 3 above). When mentioning
possible sizes of blocking sets we will always tacitly assume that
they are minimal, so deleting a point destroys the blocking
property.

A very interesting and probably feasible problem is to
characterize the sets that give equality in the bound for
$PG(2,p)$. For all (odd) $p$ the graph of the function
\[
f\,:\,x\mapsto x^{(p+1)/2},
\]
(the first example) together with it's $(p+3)/2$ directions is an
example, and it is the essentially unique one of R\'edei type
(this was proved already in 1981 by Lov\'asz and Schrijver
\cite{LS81} who also gave an elementary proof of R\'edei's result
for the case that $q=p$ is prime). Only two examples (of size
$3(p+1)/2$) are known that are not of R\'edei type, one (with 12
points) in the plane of order 7, it looks like a dual affine plane
of order 3. The other (with 21 points) in the plane of order 13
was only found last year by Blokhuis, Brouwer and Wilbrink. Both
are unique \cite{BloBrWi}. In the same paper it is shown that no
other examples exist in planes of (prime) order less than 37, and
it is extremely unlikely that this is different later on.

Motivated by these results we call blocking sets of size
$<{3\over2}(q+1)$ {\em small}, so small blocking sets only exist
in planes of non-prime order. The structure of small blocking sets
is restricted by the following theorem of Sz\H{o}nyi
\cite{Szonyi}:
\begin{theorem}{\rm [Sz\H{o}nyi, 1997]}
Let $B$ be a (minimal) small blocking set in $PG(2,q)$, where $q$
is a power of the prime $p$. Then $|B\cap l|=1$ mod $p$ for every
line $l$.
\end{theorem}

For a long time I was convinced, and even conjectured that small
blocking sets were necessarily of R\'edei type, but this turned
out to be false. Nice examples of small non-R\'edei type blocking
sets were found by Polito and Polverino \cite{PP}.

The basic idea is very simple. Consider $PG(2,q^s)$. By definition
its points and lines are the 1- and 2-dimensional subspaces of
$V=V(3,q^s)$, a 3 dimensional vector space over $GF(q^s)$. When we
consider $V'$ which is just $V$ as $3s$-dimensional over $GF(q)$
then points and lines correspond to certain $s$-, and
$2s$-dimensional subspaces of $V'$. Now let $W$ be any
$s+1$-dimensional subspace of $V'$. Let $B(W)$ be the collection
of points in $PG(2,q^s)$ for which the corresponding $s$-space in
$V'$ intersects $W$ non-trivially. One readily checks that $B(W)$
is a blocking set (of size at most $(q^{s+1}-1)/(q-1)$), because
in the $3s$ dimensional vector space $V'$ an $(s+1)$-space and a
$(2s)$-space must intersect in at least a 1-space. Polito and
Polverino give examples that are not of R\'edei type in all planes
$PG(2,p^n)$, $n>3$. The examples of small blocking sets of R\'edei
type also fall under this more general construction, by the
direction theorem.

Next we consider multiple blocking sets. $B$ is called a $t$-fold
blocking set if every line intersects $B$ in at least $t$ points.
If $q$ is a square, then $PG(2,q)$ can be partitioned into Baer
subplanes, and taking $t$ of them produces a set with the property
that every line intersects it in either $t$ or $t+\sqrt q$ points
(this makes it a two-intersection set). Again using the theory of
lacunary polynomials it can be shown that for small $t$ these are
the minimal examples \cite{BlStSz}. Our knowledge on the structure
of (relatively) small multiple blocking sets is summarized in the
following
\begin{theorem}{\rm [Blokhuis, Storme, Sz\H{o}nyi, 1998]}
Let $B$ be a $t$-fold blocking set in $PG(2,q)$ of size
$t(q+1)+c$. Let $c_2=c_3=2^{-1/3}$ and $c_p=1$ for $p>3$.
\begin{enumerate}
\item If $q=p^{2d+1}$ and $t<q/2-c_pq^{2/3}/2$ then $c\ge c_pq^{2/3}$.
\item If $4<q$ is a square, $t<q^{1/4}/2$ and $c<c_p q^{2/3}$, then
$c\ge t\sqrt{q}$ and $B$ contains the union of $t$ disjoint Baer
subplanes.
\item If $q=p^2$ and $s<q^{1/4}/2$ and $c<p\lceil {1\over4}+\sqrt{(p+1)/2}
\rceil$, then $c\ge t\sqrt{q}$ and $B$ contains the union of $t$
disjoint Baer subplanes.
\end{enumerate}
\end{theorem}

What it essentially says is that for $t<q^{1/6}$ a $t$-fold
blocking set has at least the size of $t$ disjoint Baer subplanes,
and equality implies that it is just that. In the special case
that $q$ is the square of a prime, then the same is true for
$t<q^{1/4}/2$.

This result appears to be rather sharp in the following sense: In
\cite{BaBlLa} Ball, Blokhuis and Lavrauw construct a
two-intersection set with the same parameters as, but different
from the union of $q^{1/4}+1$ disjoint Baer subplanes, so the
above characterization does no longer apply if $t>q^{1/4}$.

The construction is based on the Polito-Polverino idea. So the
plane is $PG(2,q^s)$, and $V=V(3,q^s)$ and $V'=V(3s,q)$. Now we
take for $W$ an $s+2$-dimensional subspace of $V'$. If this has
the additional property that intersections with the
$s$-dimensional subspaces of $V'$ corresponding to projective
points are at most $1$-dimensional, then the corresponding set
$B(W)$ is a $(q+1)$-fold blocking set. To see this note that if
$L$ is a $2s$-dimensional subspace of $V'$ corresponding to a
line, then $W\cap L$ is at least 2-dimensional, but by assumption
$W$ intersects $s$-spaces corresponding to points in at most 1
dimension, so it has to intersect at least $q+1$ of them.

The question of whether it is possible to find such subspaces lead
to the notion of {\em scattered subspaces} with respect to
spreads. An $s$-spread in a vector space $V$ is a collection of
$s$-dimensional subspaces partioning the nonzero vectors of $V$.
In order for $V$ to admit an $s$-spread it is necessary and
sufficient that its dimension is a multiple of $s$. So in the
above example the $s$-spaces in $V'$ corresponding to points of
$PG(2,q^s)$ define an $s$-spread. Given a vector space $V$
together with an $s$-spread $S$ we say that the subspace $W$ is
scattered by $S$ if $W$ intersect each spread element in an at
most 1-dimensional subspace. A natural question is what the
maximal dimension is of a scattered subspace. Results on this
question and related problems can be found in the thesis of
Lavrauw \cite{Lavrauw}.

A detailed survey of the many recent results on blocking sets and
multiple blocking sets is contained in the paper by Hirschfeld and
Storme \cite{HiSt2}. Blocking sets of projective planes can also
be considered as a special case of the more general concept of
covers in hypergraphs, extensively treated in the excellent (but
not too recent) survey by F\"uredi \cite{Furedi}.

To conclude this section let me mention two attractive problems
(on which no progress has been made in the last 10 years).

The first concerns double blocking sets in $PG(2,p)$. A lower
bound due to Blokhuis and Ball gives $|B|\ge 5(p+1)/2$. A trivial
example is formed by the union of three lines, of size $3p$. Could
it be that this is the  minimal size? It is true for $p=2,3,5,7$,
but it might already be false for $p=11$.

The second question has repeatedly been asked to me by Paul
Erd\H{o}s. Is there a universal constant $c$ (10 say), such that
in any plane (or any $PG(2,p)$ say) there is a blocking set with
at most $c$ points on every line. In all of the known examples
there are some lines with many points of the blocking set. Results
by Ughi show that it does not work to use for $B$ the union of a
small set of algebraic curves of bounded degree \cite{Ughi}, using
for instance a union of conics one obtains blocking sets with
$c\log(q)$ points on a line. On the other hand, in $PG(2,p^n)$
there is a blocking set with at most $p+1$ point on every line.

\section*{3. Lacunary polynomials} \addsec

\vskip-5mm \hspace{5mm}

We now turn to the main tool in the recent investigations on
(multiple) blocking sets.

Let $f(X)\in GF(q)[X]$ be fully reducible, in other words, $f(X)$
factors into linear factors over $GF(q)$. In \cite{Redei} R\'edei
investigates the case $f(X)=X^q+g(X)$ with deg$(g)<q-1$, and calls
the polynomial lacunary. The problem is to characterize those $f$
where the degree of $g$ is small. As an easy example we prove:

\begin{theorem}{\rm [R\'edei, 1970]}
Let $f(X)=X^q+g(X)$ be fully reducible in $GF(q)[X]$, where
$q=p^n$ is prime. Then either $f(X)\in GF(q)[X^p]$, or $g(X)=-X$
or deg$(g)\ge(q+1)/2$.
\end{theorem}

\begin{Proof}
Write $f=s.r$, where $s$ has the same zeroes as $f$, but with
multiplicity one, and $r$ consists of the remaining factors. Then
$s\,|\,X^q-X$, as well as $f$, so $s\,|\,X+g$, and $r\,|\,f'=g'$.
Hence $f=s.r\,|\,(X+g)g'$ so either $(X+g)g'=0$ or
$\deg(g)+\deg(g')\le \deg(f)=q$.
\end{Proof}

If $q=p$ is prime, then $g'=0$ together with deg$(g)<p$ imply that
$g$ is constant.

Much of R\'edei's book is devoted to the classification of those
$f$ with deg$(g)=(q+1)/2$. For us the case $g'=0$ (and hence $f\in
GF(q)[X^p]$ is more interesting however, also for our applications
we need to consider polynomials of the form  $f(X)=X^qg(X)+h(X)$,
where both $g$ and $h$ have degree less than $q$.

The following theorem summarizes what we know in this case
\cite{BlStSz}:

\begin{theorem}{\rm [Blokhuis, Storme, Sz\H{o}nyi, 1998]}
Let  $f\in GF(q)[X]$, $q=p^n$, $p$ prime,
 be fully reducible, $f(X)=X^qg(X)+h(X)$, where
$(g,h)=1$. Let $k=\max(\deg(g),\deg(h)<q$. Let $e$ be maximal such
that $f$ is a $p^e$-th power (so $f\in GF(q)[X^{p^e}]$). Then we
have one of the following possibilities:
\begin{enumerate}
\item
$e=n$ and $k=0$;
\item
$e\ge 2n/3$ and $k\ge p^e$;
\item
$2n/3>e>n/2$ and $k\ge p^{n-e/2}-{3\over2}p^{n-e}$;
\item
$e=n/2$ and $k=p^e$ and $f(X)=a\Tr (bX+c)+d$ or $f(X)=a\Norm
(bX+c)+d$ for suitable constants $a,b,c,d$. Here $\Tr$ and $\Norm$
respectively denote the trace and norm function from $\F_q$ to
$\F_{\sqrt q}$;
\item
$e=n/2$ and $k\ge p^e\left\lceil {1\over
4}+\sqrt{(p^e+1)/2}\right\rceil$;
\item
$n/2>e>n/3$ and $k\ge p^{(n+e)/2}-p^{n-e}-p^e/2$, or if $3e=n+1$
and $p\leq 3$, then $k\ge p^e(p^e+1)/ 2$;
\item
$n/3\ge e>0$ and $k\ge p^e\lceil (p^{n-e}+1)/(p^e+1) \rceil$;
\item
$e=0$ and $k\ge (q+1)/2$;
\item
$e=0$, $k=1$ and $f(X)=a(X^q-X)$.
\end{enumerate}
\end{theorem}

It would be very pleasant to have stronger information in the case
$n/3<e<n/2$, this would have very useful applications.

\section*{4. The connection} \addsec

\vskip-5mm \hspace{5mm}

In this section we will illustrate the connection between the
direction problem, (multiple) blocking sets, and lacunary
polynomials.

Let $f\,:\,GF(q)\to GF(q)$ be any map, and let $D_f$ be its set of
directions. Consider the auxiliary (R\'edei) polynomial

\[
R(X,Y)=\prod_{w\in GF(q)} (X-wY+f(w)),
\]
introduced by R\'edei. Let $y$ and $v\ne w\in GF(q)$. Then
$vy-f(v)=wy-f(w)$ if and only if $(f(v)-f(w))/(v-w)=y$. It follows
that for $y\not\in D_f$ the map $w\mapsto wy-f(w)$ is a bijection,
and hence $R(X,y)=\prod_{z\in GF(q)}(X-z)=X^q-X$. Write
\[
R(X,Y)=X^q+r_1(Y)X^{q-1}+\cdots+r_{q-1}(Y)X+r_q(Y),
\]
where $r_i$ is a polynomial of degree $<i$ in $Y$ (with the
exception of $r_{q-1}$ of degree $q-1$: it is clear that $r_i$ has
degree at most $i$, but the coefficient of $Y^i$ is the $i$-th
elementary symmetric polynomial in the elements of $GF(q)$, so
this is 0 for $0<i\ne q-1$). If $y\not\in D_f$ and $i\ne q-1$ then
$r_i(y)=0$. So $r_i$ is identically zero for $i\le|D_f|$. As a
consequence we have for $y\in D_f$ that $R(X,y)=X^q+g(X)$ for some
$g$ depending on $y$ with deg$(g)\le q-|D_f|-1$. So the R\'edei
polynomial when specialized for $Y=y \in D_f$ is a lacunary
polynomial and information on deg$(g)$ gives results for $|D_f|$.

The extension of the R\'edei polynomial from the graph of a
function to point sets in general can be illustrated best in the
case of an ordinary blocking set $B$ of $PG(2,q)$. We may
coordinatize our plane in such a way that the line at infinity
becomes a tangent, containing the point $(1:0:0)$ of the blocking
set.

Let $|B|=q+1+d$ then the remaining $q+d$ points have certain
affine coordinates $(a_i,b_i)$ and the R\'edei polynomial
associated to $B$ (in this position) can be defined as:
\[
R[X,Y]:=\prod(X-a_iY+b_i)=X^{q+d}+r_1(Y)X^{q+d-1}+\dots+r_{q+d}(Y).
\]
For $y,c\in GF(q)$ consider the line $\{(U,V):V=yU+c\}$. It
contains an affine point $(a,b)$ of the blocking set: $c=b-ay$.
Hence $R[X,y]$ is divisible by $(X-c)$ for all $c\in GF(q)$, in
other words $X^q-X$ divides $R[X,y]$. It follows as before that
$r_i$ is identically zero for $i=d+1,\dots,q-2$. As a consequence
we obtain the lacunary polynomial
\[
f(X)=\prod(X-a_i)=X^q g(X)+h(X),
\]
with $g$ of degree $d$ and $h$ of degree at most $d+1$, and
information on $f$ translates to information on $B$. For the case
that $q=p$ is prime it not only gives that $|B|\ge 3(p+1)/2$, but
in case of equality it also gives that each point of the blocking
set is on exactly $(p-1)/2$ tangents, and by classifying the
possible polynomials it gives all possibilities for the
configuration of the tangents through a particular point. For
small $p$ (at most 37) the number of possibilities is sufficiently
small to be handled by a computer, and to prove uniqueness of the
minimal example (for $p\ne 7, 13$).

\bibliographystyle{plain}

\label{lastpage}

\end{document}